\documentclass[11pt]{article}

\RequirePackage{amsthm,amsmath,amsfonts,amssymb}
\RequirePackage[authoryear]{natbib}

\usepackage{todonotes,soul}

\usepackage[T1]{fontenc}
\usepackage[utf8]{inputenc}
\usepackage{lmodern}
\usepackage{textcomp}

\usepackage[english]{babel}

\usepackage{enumitem}

\usepackage{color}

\usepackage[margin=2.8cm]{geometry}

\usepackage{url}
\usepackage[unicode,colorlinks,pdfusetitle,final]{hyperref}

\usepackage[final]{microtype}

\theoremstyle{plain}
\newtheorem{theorem}{Theorem}

\newtheorem{proposition}[theorem]{Proposition}

\theoremstyle{definition}

\theoremstyle{remark}
\newtheorem{remark}[theorem]{Remark}

\newcommand{\field}[1][K]{\mathbb{#1}}
\newcommand{\NN}{\field[N]}

\newcommand{\RR}{\field[R]}

\newcommand{\Bc}{\mathcal{B}}
\newcommand{\Cc}{\mathcal{C}}
\newcommand{\Beta}{{\sf Beta}}
\newcommand{\Uniform}{{\sf Uniform}}

\newcommand{\prob}[1]{\mathbb{P}(#1)}
\newcommand{\pFq}[5]{{}_{#1}F_{#2}\left(\begin{array}{c} #3 \\ #4 \end{array}; #5\right)}

\definecolor{hypercolor}{rgb}{0,0.2,0.7}
\hypersetup{
  linkcolor=hypercolor,
  urlcolor=hypercolor,
  citecolor=hypercolor
}

\newcommand{\oeis}[1]{\href{http://oeis.org/#1}{\texttt{#1}}}

\newcommand{\KdF}[9]{F^{#1}_{#2}\left(\begin{array}{c} #3: #4 \mathbin{;} #5 \\ #6: #7\mathbin{;} #8 \end{array}; #9\right)}

\begin{document}

\title{Integrals of incomplete beta functions, with applications to order statistics, random walks and string enumeration}

\author{%
Stephen B. Connor\footnote{\tt stephen.connor@york.ac.uk}\; and Christopher J. Fewster\footnote{\tt chris.fewster@york.ac.uk}  \\
\small \emph{Department of Mathematics, University of York, Heslington, York, YO10 5DD, UK}
}
\date{\today} 
\maketitle

\begin{abstract}
We study the probability that one beta-distributed random variable exceeds the maximum of two others,
allowing all three to have general parameters. This amounts to studying Euler transforms of products of two incomplete beta functions. We provide a closed form for the general problem in terms of 
Kamp\'e de F\'eriet functions and a variety of simpler closed forms in special cases. The results are applied to derive the moments of the maximum of two independent beta-distributed random variables and to find inner products of incomplete beta functions. Restricted to positive integer parameters, our results are applied to determine an expected exit time for a conditioned random walk and also to a combinatorial problem of enumerating strings comprised of three different letters, subject to constraints. 


\smallskip

\end{abstract}

\section{Introduction} \label{sec:intro}
In this note we investigate and document various integrals of incomplete beta functions,
of interest in their own right, which have applications to probability theory and combinatorics. Among other things, we will provide formulae for: (a) all moments of the maximum of two independent random variables following beta distributions with arbitrary parameters, (b) inner products of incomplete beta functions, and (c) solutions to a combinatorial problem involving the enumeration of strings of certain type.

For parameters $\mu,\mu'>0$, we write $X_{\mu,\mu'}\sim\Beta(\mu,\mu')$ for the random variable with density proportional to $s^{\mu-1}(1-s)^{\mu'-1}$ for $s\in[0,1]$. The cumulative distribution function of $X_{\mu,\mu'}$ is given for $t\in[0,1]$ by the regularized incomplete beta function
\[ 
I(\mu,\mu';t) =  \prob{X_{\mu,\mu'} \le t}  = \frac{1}{B(\mu,\mu')}\int_0^t s^{\mu-1}(1-s)^{\mu'-1}\, ds\,,
\]
where $B(\mu,\mu')$ is the beta function
\[ 
B(\mu,\mu') = \frac{\Gamma(\mu) \Gamma(\mu')}{\Gamma(\mu+\mu')}  = \int_0^1 s^{\mu-1}(1-s)^{\mu'-1}\, ds \,.
\]

We are interested in calculating integrals of the following type:
\begin{equation}\label{eq:B_int_defn}
\Bc(\lambda,\lambda',\mu,\mu',\nu,\nu'):=\frac{1}{B(\lambda,\lambda')} \int_0^1  t^{\lambda-1}(1-t)^{\lambda'-1} I(\mu,\mu';t) I(\nu,\nu';t)\,dt    \,,
\end{equation}
and will provide a closed form for the general case, in terms of Kamp\'e de F\'eriet functions~\citep{KampedeFeriet1921a,KampedeFeriet1921b}, and a variety of special cases in terms of
hypergeometric and elementary functions.
When the six parameters are all positive reals, the integrals converge and have the probabilistic interpretation
\begin{equation}\label{eq:prob_interpretation}
\Bc(\lambda,\lambda',\mu,\mu',\nu,\nu') = \prob{X_{\lambda,\lambda'}>\max\{Y_{\mu,\mu'},Z_{\nu,\nu'}\}} \,,
\end{equation}
where $X_{\lambda,\lambda'},Y_{\mu,\mu'},Z_{\nu,\nu'}$ are independent beta random variables, as is clear on noting that $I(\mu,\mu';t) I(\nu,\nu';t)$ is the cumulative distribution function for $\max\{Y_{\mu,\mu'},Z_{\nu,\nu'}\}$.  In consequence the following identities are obvious:
\begin{equation}\label{eq:swap}
\Bc(\lambda,\lambda',\mu,\mu',\nu,\nu') = \Bc(\lambda,\lambda',\nu,\nu',\mu,\mu')
\end{equation}
and
\begin{equation}\label{eq:cyclicity}
\Bc(\lambda,\lambda',\mu,\mu',\nu,\nu') + \Bc(\mu,\mu',\nu,\nu',\lambda,\lambda') + \Bc(\nu,\nu',\lambda,\lambda',\mu,\mu')=1\,.
\end{equation}
Furthermore, setting $\lambda' = 1$ yields
\begin{align}\label{eqn:Bc-moments}
\Bc(\lambda,1, \mu,\mu',\nu,\nu') &= \lambda\int_0^1 t^{\lambda-1}I(\mu,\mu';t) I(\nu,\nu';t)\,dt  \notag \\
&= 1- \int_0^1 t^\lambda\, \frac{d}{dt}\left(I(\mu,\mu';t)I(\nu,\nu';t)\right)\, dt \notag \\
&=	1 - \mathbb{E}[\max\{Y_{\mu,\mu'},Z_{\nu,\nu'}\}^\lambda] \,,
\end{align}
and so all moments of the random variable $\max\{Y_{\mu,\mu'},Z_{\nu,\nu'}\}$ may be deduced from the family of integrals in~\eqref{eq:B_int_defn}. Thanks to the invariance of the beta distribution under the transformation $(s,\mu)\leftrightarrow (1-s,\mu')$, $\Bc(\lambda,\lambda',\mu,\mu',\nu,\nu')$ may equivalently be expressed as
\begin{equation}\label{eq:min_not_max}
\Bc(\lambda,\lambda',\mu,\mu',\nu,\nu')=	\prob{X_{\lambda',\lambda}<\min\{Y_{\mu',\mu},Z_{\nu',\nu}\}}\,;
\end{equation}
hence our results can just as well be used to determine properties of the minimum of two independent beta random variables. 

A closed form for $\Bc(\lambda,1, \mu,\mu',\nu,\nu')$ in terms of ${}_3F_2$ hypergeometric functions at unit argument (defined in \eqref{eq:hypgeo}) is given in~\eqref{eq:lambda'=1} below and immediately yields the following expectation values.
\begin{proposition}\label{prop:moments}
	Let $Y_{\mu,\mu'}$ and $Z_{\nu,\nu'}$ be independent beta random variables as defined above. Then for any $\lambda\ge 0$ the $\lambda$'th moment of $\max\{Y_{\mu,\mu'},Z_{\nu,\nu'}\}$ is given by
	\begin{align}
	\mathbb{E}[\max\{Y_{\mu,\mu'},Z_{\nu,\nu'}\}^\lambda]&=
	\frac{B(\lambda+\mu+\nu,\mu')}{\nu B(\mu,\mu') B(\nu,\nu')}
	\,\pFq 32{1-\nu',\nu, \lambda+\mu+\nu}{\nu+1,  \mu'+\lambda+\mu+\nu}{1} \nonumber\\
	&\qquad + \frac{B(\lambda+\mu+\nu,\nu')}{\mu B(\mu,\mu') B(\nu,\nu')}
	\,\pFq 32{1-\mu',\mu, \lambda+\mu+\nu}{\mu+1,  \nu'+\lambda+\mu+\nu}{1}\,.
	\label{eq:EmaxYZlambda}
	\end{align}
\end{proposition}

The random variable $\max\{Y_{\mu,\mu'},Z_{\nu,\nu'}\}$ does not follow a standard distribution, yet particular cases of the integral in \eqref{eq:B_int_defn} arise in a number of rather diverse situations. For example, the Program Evaluation and Review Technique (PERT) is a widely-used method for analysing the time-to-completion of a project involving multiple subtasks. Independent beta random variables are commonly used to model the time taken for each individual subtask to be completed; the maximum of two such distributions is then sometimes approximated by another beta distribution with (numerically approximated) matching first and second moments \citep{sculli1985}. 

Calculating properties of the maximum of a set of independent random variables falls, of course, under the more general umbrella of determining their order statistics. In the case of beta random variables, methods are provided in \cite{thomas2008} and \cite{abdelkader2010} for computing moments of order statistics in the case where the parameters of the beta distribution(s) are all integer-valued. However, the resulting equations are either recursive~\citep{thomas2008} or quite complicated nested sums of `inclusion exclusion' type~\citep{abdelkader2010}. For example, in \cite{thomas2008} (in which it is also assumed that the random variables are all identically distributed) the expectation of the maximum of two $\Beta(m,n)$ random variables is expressed as a linear combination of the first $m+n$ moments of a single $\Beta(m,n)$. 
The work of \cite{nadarajah2008} drops the assumption of integer-valued parameters but, unlike the problem considered in the present paper, once again requires that the random variables are identically distributed. In that setting it is shown that the moments of beta order statistics may be expressed as a sum of Kamp\'e de F\'eriet functions.

By contrast, in Section~\ref{ssec:special} we shall show that in the case of two i.i.d. random variables, Proposition~\ref{prop:moments} leads to the simple result 
\begin{equation}\label{eq:B11munumunu}
\mathbb{E}[\max\{Y_{\mu,\nu},Z_{\mu,\nu}\}] = 1-\Bc(1,1,\mu,\nu,\mu,\nu)=
\frac{\mu}{\mu+\nu} + \frac{2B(2\mu,2\nu)}{(\mu+\nu)B(\mu,\nu)^2}
\end{equation}
for arbitrary real $\mu,\nu>0$ (see~\eqref{eq:B11mnmn}); we note that the first term here is equal to $\mathbb{E}[Y_{\mu,\nu}]$, and the second to $\tfrac12 \mathbb{E}[|Y_{\mu,\nu}-Z_{\mu,\nu}|]$. Similarly, when $Y$ and $Z$ both follow symmetric beta distributions we have the following formula:
\begin{equation}\label{eq:B11mumununu}
\mathbb{E}[\max\{Y_{\mu,\mu},Z_{\nu,\nu}\}] 
= \frac{1}{2} +\frac{(\mu+\nu)B(\mu+\nu,\mu+\nu+1)}{4\mu\nu B(\mu,\mu+1)B(\nu,\nu+1)} \,.
\end{equation}

\cite{cordeiro2011} introduce a number of beta generalized distributions and study their properties. This family includes the so-called \textit{beta beta distribution} ${\sf BB}(a,b,\mu,\mu')$, with density  
\[
\frac{t^{\mu-1}(1-t)^{\mu' -1}}{B(\mu,\mu')B(a,b)}I(\mu, \mu';t)^{a-1}(1-I(\mu, \mu';t))^{b-1}\,, \quad t\in[0,1].
\]
When $a=b=2$ this density is similar to the integrand in \eqref{eq:B_int_defn}, but requires the three beta (or incomplete beta) functions to have identical parameters. It is shown by  \cite{cordeiro2011} that the moments of the beta beta distribution may be expressed as infinite sums of generalized Kamp\'e de F\'eriet functions. The particular case where $a=b=2$ may be rewritten (and, in some cases, greatly simplified) using our results. Indeed, if $W\sim {\sf BB}(2,2,\mu,\mu')$ then for any $\lambda\ge 0$, 
\begin{align*}
\frac16 \mathbb{E}[W^\lambda] &=  	\int_0^1 \frac{t^{\lambda+\mu-1}(1-t)^{\mu' -1}}{B(\mu,\mu')}I(\mu, \mu';t)(1-I(\mu, \mu';t)) \,dt\\
&= \prob{Y_{\mu,\mu'}< X_{\lambda+\mu,\mu'} < Z_{\mu,\mu'}} \\	
&= \prob{Y_{\mu,\mu'}< X_{\lambda+\mu,\mu'}} - \Bc(\lambda+\mu,\mu',\mu,\mu',\mu,\mu') 	\,,
\end{align*}
and these two terms may be computed using results derived below (Remark~\ref{rem:simpler_problem} and Proposition~\ref{prop:general_result}, respectively).

\medskip
The integrals $\Bc(\lambda,\lambda',\mu,\mu',\nu,\nu')$ are of course interesting in their own right as properties of special functions. For general positive real parameters, we will first of all give a formula in \eqref{eq:4F3series} exhibiting $\Bc(\lambda,\lambda',\mu,\mu',\nu,\nu')$ as a series of terms involving ${}_4F_3$ hypergeometric functions at unity that terminates if $\mu'\in\NN$ (by symmetry~\eqref{eq:swap}, there is a similar formula terminating if $\nu'\in\NN$). We then show that a little more work yields the following slightly simpler general formula involving ${}_3F_2$ hypergeometric functions.

\begin{proposition}\label{prop:general_result}
	For any set of positive real-valued arguments, the integral $\Bc(\lambda,\lambda',\mu,\mu',\nu,\nu')$ is given by
	\[
	\Bc(\lambda,\lambda',\mu,\mu',\nu,\nu')=1 - 
	\frac{\Cc(\lambda,\lambda',\mu,\mu',\nu,\nu')+\Cc(\lambda,\lambda',\nu,\nu',\mu,\mu')}{B(\lambda,\lambda')B(\mu,\mu')B(\nu,\nu')}
	\]
	where
	\begin{equation}\label{eq:Cseries}
	\Cc(\lambda,\lambda',\mu,\mu',\nu,\nu') :=
	\sum_{k=0}^\infty \frac{(1-\lambda')_k\,B(\lambda+k+\mu+\nu,\mu')}{k!(\lambda+k)\nu} 
	\,\pFq 32{1-\nu',\nu, k+\lambda+\mu+\nu}{\nu+1,  k+\mu'+\lambda+\mu+\nu}{1}    
	\end{equation}
	terminates if $\lambda'\in\NN$, and for general parameters can be written in terms of a Kamp\'e de F\'eriet function:
	\[
	\Cc(\lambda,\lambda',\mu,\mu',\nu,\nu')=
	\frac{B(\lambda+\mu+\nu,\mu')}{\lambda\nu} 
	\KdF{1:2;2}{1:1;1}{\lambda+\mu+\nu}{\lambda,1-\lambda'}{\nu,1-\nu'}{\lambda+\mu+\nu+\mu'}{\lambda+1}{\nu+1}{1,1}\,. 
	\]
\end{proposition}
We will also show that there is a similar expansion that terminates if
$\lambda\in\NN$; in the light of~\eqref{eq:cyclicity} it follows that
any $\Bc(\lambda,\lambda',\mu,\mu',\nu,\nu')$ in which $\{\mu,\mu'\}$ and $\{\nu,\nu'\}$ both contain at least one natural number may be expressed in terms of a finite sum of ${}_3F_{2}$'s at unit argument. The general definition of the Kamp\'e de F\'eriet function is given in~\eqref{eq:KdFdefn}. It is easily implemented, for instance, in the high precision Python library {\tt mpmath} using the {\tt hyper2d} function; sample code for the function $\Cc$ is given below: 

\bigskip
\hrule
\begin{verbatim}
import mpmath as mp

def C(l,lpr,m,mpr,n,npr):
lmn=l+m+n
return mp.beta(lmn,mpr)/(l*n)*mp.hyper2d({'m+n':[lmn], 'm':[1-lpr,l],
'n':[1-npr,n]},{'m+n':[lmn+mpr], 'm':[l+1], 'n':[n+1]}, 1,1)	
\end{verbatim}
\hrule
\bigskip

When $\lambda=\lambda' = 1$, our results may alternatively be interpreted as the inner product of $I(\mu,\mu';t)$ and $I(\nu,\nu';t)$ in $L^2((0,1),dt)$, with Proposition~\ref{prop:moments} showing that 
\begin{align}\label{eq:IPs}
\int_0^1 I(\mu,\mu';t) I(\nu,\nu';t)\,dt &= 1 - 
\frac{B(1+\mu+\nu,\mu')}{\nu B(\mu,\mu') B(\nu,\nu')}
\,\pFq 32{1-\nu',\nu, 1+\mu+\nu}{\nu+1,  \mu'+1+\mu+\nu}{1}\notag
\\&\qquad -  \frac{B(1+\mu+\nu,\nu')}{\mu B(\mu,\mu') B(\nu,\nu')}
\,\pFq 32{1-\mu',\mu, 1+\mu+\nu}{\mu+1,  \nu'+1+\mu+\nu}{1}\,.
\end{align}
Inner products of this type have useful applications. For example, 
our original motivation to study this problem arose from the second author's desire to calculate the $L^2$-norm of the function $p_\mu(t)=I(\mu,\mu;t)^2$ \citep{fewster2019}. In that application -- which concerned singularity theorems in general relativity generalising those of \cite{Hawking:1966} -- incomplete beta functions enter as mollified step functions. If $\mu=m\in\NN$, one has $p_m(t)=O(t^m)$ as $t\to 0^+$ and $p_m(t)=1-O((1-t)^m)$ as $t\to 1^-$, with $p_m$ monotone increasing on $[0,1]$. Extending by zero to the left and unity to the right, one obtains a $C^{m-1}$ mollified step function; concatenating with a reflected version of itself, one obtains a bump function belonging to the Sobolev space $W_0^{m}(\RR)$. The application required $\|p_m\|$, $\|p_m'\|$ and $\|p_m^{(m)}\|$, where the norms are those of $L^2((0,1),dt)$. The formulae
\[ 
\|p'_m\|^2  = \frac{B(2m-1,2m-1)}{B(m,m)^2}, \qquad
\|p_m^{(m)}\|^2=   \frac{(2m-2)!(2m-1)!}{(m-1)!^2} 
\]
are straightforward~\cite[Appendix]{fewster2019}, while the formula 
\[ 
\|p_m\|^2 =\Bc(1,1,m,m,m,m)= \frac{1}{2} - \frac{(2m)!^4}{4(4m)!m!^4}
\]
(stated but not proved in~\cite{fewster2019}) is obtained by the methods of the present note.
Indeed, writing~\eqref{eq:B11mumununu} for integer parameters one has the inner product
\begin{equation}\label{eq:B11mmnn}
\langle p_m , p_n \rangle=\Bc(1,1,m,m,n,n)  = \frac12 - \frac14 \frac{\binom{m+n}{m}^2}{ \binom{2m+2n}{2m} } \qquad (m,n\in\NN)\,,
\end{equation}
from which the formula for $\|p_m\|^2$ given above follows immediately. 

As will be shown in Section~\ref{sec:general}, $\Bc(\lambda,\lambda',\mu,\mu',\nu,\nu')$ may be recouched using generalised hypergeometric functions. A recent paper of \cite{conway2017} deals with some related indefinite integration formulae, but unfortunately these do not include our integral of interest.
The methods we will use are essentially standard, the main tool being the Euler transform of hypergeometric functions, but some of the special cases are evaluated using `contiguous' generalisations of the classical Whipple and Dixon identities for hypergeometric functions derived by \cite{Lavoie:1994} and \cite{Lavoie:1996}. In Section~\ref{sec:path_counting} we specialise to integer parameters and give two applications of our results. The first concerns exit times in a conditioned random walk. Specifically, we consider a random walk on a rectangular lattice, starting from the origin, taking northward or eastward steps with equal probability. Conditioned on the path passing through some point $(M,N)$, we show that the expected exit time from the rectangle 
with corners $(0,0)$ and $(M-m,N-n)$ is $(M+N+1)\Bc(1,1,m,M-m+1,n,N-n+1)$, where $1\le m\le M$, $1\le n\le N$. Second, we show that $\Bc(\ell,\ell',m,m',n,n')$ determines the number of strings that can be constructed using three characters {\tt X, Y, Z}, repeated $\ell+\ell'-1$, $m+m'-1$ and $n+n'-1$ times respectively, subject to the constraint that the $\ell$'th {\tt X} appears after at least $m$ {\tt Y}'s and $n$ {\tt Z}'s. Thus our results provide closed form solutions to this combinatorial problem. Various integer sequences that emerge from special cases of this analysis are presently unknown to the On-Line Encyclopedia of Integer Sequences,\footnote{\href{https://oeis.org/}{https://oeis.org/}} indicating the novelty of our results.

\section{Main results}\label{sec:general}

\subsection{General formulae}
We begin by recalling some basic facts concerning generalised hypergeometric functions, defined as usual by
\begin{equation}\label{eq:hypgeo}
\pFq{p}{q}{a_1,\dots,a_p}{b_1,\dots,b_q}{t} = \sum_{k=0}^\infty \frac{(a_1)_k \dots (a_p)_k}{(b_1)_k \dots (b_q)_k} \, \frac{t^k}{k!}\,,
\end{equation}
where the Pochhammer symbol $(a)_k = a(a+1)\cdots(a+k-1)$ denotes the $k^{th}$ rising factorial of $a$. The Euler transform (\citeauthor[Eq.16.5.2]{DLMF}) states that
\begin{equation}\label{eq:Euler_transform}
\int_0^1  u^{c-1}(1-u)^{d-c-1}\, \pFq{p}{q}{a_1,\dots,a_p}{b_1,\dots,b_q}{tu}\,du  
=  \frac{\Gamma(c)\Gamma(d-c)}{\Gamma(d)}\,\pFq{p+1}{q+1}{a_1,\dots,a_p,c}{b_1,\dots,b_q,d}{t}\,.
\end{equation}
Starting with ${}_1F_0(-a;t)=(1-t)^a$ the Euler transform immediately yields one of the hypergeometric representations of the regularized incomplete beta function:
\begin{align}
I(\mu,\mu';t) &= \frac{1}{B(\mu,\mu')} \int_0^t s^{\mu-1}(1-s)^{\mu'-1}\,ds = \frac{t^\mu}{B(\mu,\mu')} \int_0^1 u^{\mu-1}
\,{}_1F_0(1-\mu';tu)\,du \nonumber \\
&= \frac{t^\mu}{\mu B(\mu,\mu')} \,\pFq 21{1-\mu', \mu}{\mu+1}{t} \,. \label{eq:I_as_pFq}
\end{align}
Accordingly, our integral $\Bc(\lambda,\lambda',\mu,\mu',\nu,\nu')$ may be expressed in terms of an Euler transform of a product of hypergeometric functions
\begin{align}\label{eq:ETFF}
\Bc(\lambda,\lambda',\mu,\mu',\nu,\nu') &= 
\frac{1}{\mu\nu B(\lambda,\lambda')B(\mu,\mu')B(\nu,\nu')} \nonumber \\
&\qquad\times
\int_0^1 t^{\lambda+\mu+\nu-1}(1-t)^{\lambda'-1} \pFq 21{1-\mu', \mu}{\mu+1}{t}
\pFq 21{1-\nu', \nu}{\nu+1}{t}\,dt\,.
\end{align} 
A product of ${}_2F_{1}$ hypergeometric functions of $t$ can be represented as a series in $t$ with coefficients involving values of ${}_4F_3$ hypergeometric functions evaluated at $1$ \citep[4.3(14)]{Erdelyi_etal_vol1:1953}, giving, in this case,  
\[
\pFq 21{1-\mu',~ \mu}{\mu+1}{t}
\pFq 21{1-\nu',~ \nu}{\nu+1}{t} = 
\sum_{n=0}^\infty \frac{\mu(1-\mu')_n t^n}{(\mu+n)n!} 
\pFq 43{1-\nu',~\nu,~-\mu-n,~-n}{1+\nu,~\mu'-n,~1-\mu-n}{1}\,.
\]
The two ${}_2F_{1}$ series converge absolutely on the unit circle (see~\citeauthor[16.2(iii)]{DLMF}), and so their product in turn has an absolutely convergent Taylor series on the unit circle. Fubini's theorem thus permits the interchange of sum and integral once the series above has been substituted into \eqref{eq:ETFF}, yielding the general formula
\begin{align}\label{eq:4F3series}
\Bc(\lambda,\lambda',\mu,\mu',\nu,\nu') &= 
\frac{1}{\nu B(\lambda,\lambda')B(\mu,\mu')B(\nu,\nu')}
\sum_{n=0}^\infty \frac{(1-\mu')_n}{(\mu+n)n!}B(\lambda+\mu+\nu+n,\lambda')\nonumber \\
&\qquad\times
\pFq 43{1-\nu',~\nu,~-\mu-n,~-n}{1+\nu,~\mu'-n,~1-\mu-n}{1}\,,
\end{align}
which terminates if $\mu'\in\NN$. 
\medskip

Although its derivation was very direct, the final expression~\eqref{eq:4F3series} is somewhat unwieldy.
The simpler formula given in Proposition~\ref{prop:general_result} 
is better adapted to closed form evaluation in special cases and can be obtained as follows.
First consider the case when $\nu'=1$, with the other five parameters arbitrary. Since $I(\nu,1;t) = t^\nu$, using \eqref{eq:I_as_pFq} we see that 
\begin{align}
\Bc(\lambda,\lambda',\mu,\mu',\nu,1) &= 
\frac{1}{B(\lambda,\lambda')}\int_0^1 t^{\lambda+\nu-1}(1-t)^{\lambda'-1}I(\mu,\mu';t)\,dt \nonumber\\
&=\frac{1}{\mu B(\mu,\mu')B(\lambda,\lambda')} \int_0^1 t^{\lambda+\mu+\nu-1}(1-t)^{\lambda'-1}
\pFq 21{1-\mu', \mu}{\mu+1}{t}\,dt \nonumber\\
&= \frac{B(\lambda+\mu+\nu,\lambda')}{\mu B(\lambda,\lambda')B(\mu,\mu')}
\,\pFq 32{1-\mu', \mu, \lambda+\mu+\nu}{\mu+1, \lambda'+\lambda+\mu+\nu}{1} \,, \label{eq:nu'=1}
\end{align}
with the final equality following from a second application of the Euler transform. Using~\eqref{eq:cyclicity} and~\eqref{eq:swap} one has
\begin{align}
\Bc(\lambda,1,\mu,\mu',\nu,\nu')&= 1 - \Bc(\mu,\mu',\nu,\nu',\lambda,1) - \Bc(\nu,\nu',\mu,\mu',\lambda,1) \nonumber \\
&=1 - \frac{B(\lambda+\mu+\nu,\mu')}{\nu B(\mu,\mu') B(\nu,\nu')}
\,\pFq 32{1-\nu',\nu, \lambda+\mu+\nu}{\nu+1,  \mu'+\lambda+\mu+\nu}{1} \nonumber\\
&\qquad - \frac{B(\lambda+\mu+\nu,\nu')}{\mu B(\mu,\mu') B(\nu,\nu')}
\,\pFq 32{1-\mu',\mu, \lambda+\mu+\nu}{\mu+1,  \nu'+\lambda+\mu+\nu}{1}\,. \label{eq:lambda'=1}
\end{align} 
In passing, we have proved Proposition~\ref{prop:moments} as a trivial consequence of~\eqref{eq:lambda'=1} and~\eqref{eqn:Bc-moments}.

\medskip
Now consider the more general case of the integral in \eqref{eq:B_int_defn}, in which all parameters are real and positive.  Using the binomial theorem,
\[
\Bc(\lambda,\lambda',\mu,\mu',\nu,\nu')=\frac{1}{B(\lambda,\lambda')} \int_0^1 
\sum_{k=0}^\infty \frac{(1-\lambda')_k}{k!} t^{\lambda+k-1} I(\mu,\mu';t) I(\nu,\nu';t)\,dt   \,. 
\]
The summands are integrable functions which, for sufficiently large $k$, all have a common sign; we may therefore interchange sum and integral by Tonelli's theorem, to give
\begin{equation}\label{eq:convex_combination}
\Bc(\lambda,\lambda',\mu,\mu',\nu,\nu')=\frac{1}{B(\lambda,\lambda')}
\sum_{k=0}^\infty \frac{(1-\lambda')_k}{k!(\lambda+k)} \Bc(\lambda+k,1,\mu,\mu',\nu,\nu')\,.
\end{equation}
Using \eqref{eq:lambda'=1}, together with the identity
\begin{equation}\label{eq:beta_as_series}
\sum_{k=0}^\infty \frac{(1-\lambda')_k}{k!(\lambda+k)} = B(\lambda,\lambda')\,,
\end{equation}
which arises from a binomial expansion of the beta function integral, we see that our integral of interest can be expressed as a linear combination of ${}_3F_2$ hypergeometric functions:
\[
\Bc(\lambda,\lambda',\mu,\mu',\nu,\nu')=1 - 
\frac{\Cc(\lambda,\lambda',\mu,\mu',\nu,\nu')+\Cc(\lambda,\lambda',\nu,\nu',\mu,\mu')}{B(\lambda,\lambda')B(\mu,\mu')B(\nu,\nu')}
\]
where $\Cc(\lambda,\lambda',\mu,\mu',\nu,\nu')$ is as defined by the series given in \eqref{eq:Cseries},
which obviously terminates in the case $\lambda'\in \NN$.

For general (real positive) parameters, the series~\eqref{eq:Cseries} may be evaluated in terms of a Kamp\'e de F\'eriet function \citep{KampedeFeriet1921a,KampedeFeriet1921b}, on noting that
\begin{align*}
\frac{\lambda\nu\Cc(\lambda,\lambda',\mu,\mu',\nu,\nu')}{B(\lambda+\mu+\nu,\mu')} &= 
\sum_{k=0}^\infty \frac{(1-\lambda')_k\,(\lambda)_k\,(\lambda+\mu+\nu)_k}{k!(\lambda+1)_k\,(\lambda+\mu+\nu+\mu')_k} 
\,\pFq 32{1-\nu',\nu, k+\lambda+\mu+\nu}{\nu+1,  k+\mu'+\lambda+\mu+\nu}{1} \notag \\
&= \sum_{k=0}^\infty\sum_{j=0}^\infty \frac{(1-\lambda')_k\,(\lambda)_k\,(\lambda+\mu+\nu)_k\,(\lambda+\mu+\nu+k)_j\,(1-\nu')_j\,(\nu)_j}{(\lambda+1)_k\,(\lambda+\mu+\nu+\mu')_k\,(\lambda+\mu+\nu+\mu'+k)_j\,(\nu+1)_j\,    k! j!} \notag \\
&=\sum_{k=0}^\infty\sum_{j=0}^\infty \frac{(\lambda+\mu+\nu)_{j+k}\,(1-\lambda')_k\,(\lambda)_k\,(1-\nu')_j\,(\nu)_j}{(\lambda+\mu+\nu+\mu')_{j+k}\,(\lambda+1)_k\,(\nu+1)_j\,    k! j!}	
\end{align*}
from which
\begin{equation*}
\Cc(\lambda,\lambda',\mu,\mu',\nu,\nu') = 
\frac{B(\lambda+\mu+\nu,\mu')}{\lambda\nu} 
\KdF{1:2;2}{1:1;1}{\lambda+\mu+\nu}{\lambda,1-\lambda'}{\nu,1-\nu'}{\lambda+\mu+\nu+\mu'}{\lambda+1}{\nu+1}{1,1} 
\end{equation*}
follows by the definition of the Kamp\'e de F\'eriet function. Here, we recall that the general Kamp\'e de F\'eriet function in $2$ variables is given by
\begin{equation}
\label{eq:KdFdefn} 
\KdF{A:C;F}{B:D;G}{(a)}{(c)}{(f)}{(b)}{(d)}{(g)}{x,y}  = 
\sum_{m=0}^\infty\sum_{n=0}^\infty\frac{((a))_{m+n} ((c))_m ((f))_n}{((b))_{m+n}((d))_m ((g))_n }\, \frac{x^m y^n}{m!n!}
\,,
\end{equation}
where $(a)$, $(c)$ etc are sequences of length $A$, $C$ etc, and the notation $((c))_m=\prod_{i=1}^C (c_i)_m$ denotes the product of Pochhammer symbols. 
Thus the proof of Proposition~\ref{prop:general_result} is complete and we have evaluated the full family of integrals $\Bc(\lambda,\lambda',\mu,\mu',\nu,\nu')$ in terms of two Kamp\'e de F\'eriet functions.

\begin{remark} \label{rem:simpler_problem}
	At this point, it is useful to note the simpler problem of calculating the probability that one beta random variable exceeds another:
	\begin{align}
	\mathbb{P}(X_{\lambda,\lambda'} \ge Y_{\mu,\mu'})  &= \frac{1}{B(\lambda,\lambda')}\int_0^1 t^{\lambda-1}(1-t)^{\lambda'-1}I(\mu,\mu';t) \,dt \nonumber \\
	&=  \frac{1}{\mu B(\lambda,\lambda')B(\mu,\mu')}\int_0^1 t^{\lambda+\mu-1}(1-t)^{\lambda'-1}\,{}_2F_1 
	\left(\begin{array}{c} 1-\mu', \mu  \\ \mu+1 \end{array}; t\right) \,dt \nonumber \\
	&=  \frac{B(\lambda+\mu,\lambda')}{\mu B(\lambda,\lambda')B(\mu,\mu')}
	\,{}_3F_2 
	\left(\begin{array}{c} 1-\mu', \mu, \lambda+\mu \\ \mu+1, \lambda+\lambda'+\mu \end{array}; 1\right) \,,  \label{eq:simpler-problem}
	\end{align}
	which follows by a simple application of first \eqref{eq:I_as_pFq} and then \eqref{eq:Euler_transform}. As $\nu\to 0+$ the random variable $Z_{\nu,1}$ converges in probability to 0, and \eqref{eq:simpler-problem} agrees as expected with the corresponding limit of $\Bc(\lambda,\lambda',\mu,\mu',\nu,1)$ in \eqref{eq:nu'=1}.
\end{remark}

To conclude our general discussion, we note that \eqref{eq:min_not_max} implies that
\begin{align*}
\Bc(\lambda,\lambda',\mu,\mu',\nu,\nu') &= \prob{X_{\lambda',\lambda}<\min\{Y_{\mu',\mu},Z_{\nu',\nu}\}} \\
&= 1- \prob{X_{\lambda',\lambda}\ge Y_{\mu',\mu}} - \prob{X_{\lambda',\lambda}\ge Z_{\nu',\nu}} + \prob{X_{\lambda',\lambda}\ge\max\{Y_{\mu',\mu},Z_{\nu',\nu}\}}  \\
&= 1- \prob{X_{\lambda',\lambda}\ge Y_{\mu',\mu}} - \prob{X_{\lambda',\lambda}\ge Z_{\nu',\nu}} + \Bc(\lambda',\lambda,\mu',\mu,\nu',\nu) \,,
\end{align*} 
and, combined with Remark~\ref{rem:simpler_problem}, one may obtain from  Proposition~\ref{prop:general_result} an expression for $\Bc(\lambda,\lambda',\mu,\mu',\nu,\nu')$ involving a series of ${}_3F_{2}$'s that terminates if $\lambda\in\NN$. Combining these observations with the identity in \eqref{eq:cyclicity} shows that any $\Bc(\lambda,\lambda',\mu,\mu',\nu,\nu')$ in which both $\{\mu,\mu'\}$ and $\{\nu,\nu'\}$ contain at least one natural number may also be expressed as a finite sum of ${}_3F_{2}$'s, as claimed in Section~\ref{sec:intro}. We refrain from writing these expressions explicitly.

\subsection{Some special cases}\label{ssec:special}

We conclude this section by considering some particular sets of parameter values for which the integral $\Bc(\lambda,\lambda',\mu,\mu',\nu,\nu')$ can be expressed in terms of a small number of hypergeometric or simpler functions.\\


\noindent\textbf{Case 1:} $\lambda=\lambda'=1$, $\mu'=\mu$, $\nu'=\nu$. \\
We first consider the situation in which the random variables $Y_{\mu,\mu}$ and $Z_{\nu,\nu}$ both follow symmetric beta distributions. Here the hypergeometric functions in \eqref{eq:lambda'=1} take the form
\[
\pFq 32{a, b, c}{e, f}{1} 
\]
where $a+b=1+i+j$, $e+f=2c+1+i$, $i=-1$, $j=1$. Consequently, they may be evaluated using a contiguous generalization of Whipple's theorem proved in~\cite{Lavoie:1996}: i.e., for $a+b=1$, $e+f=2c$
\begin{align*}
\pFq 32{a, b, c}{e, f}{1} &= \frac{\Gamma(e)\Gamma(f)}{4^a(c-1)\Gamma(e-a)\Gamma(f-a)} \\
&\qquad \times\left(
\frac{\Gamma(\tfrac{1}{2}(e-a+1))\Gamma(\tfrac{1}{2}(f-a))}{\Gamma(\tfrac{1}{2}(e+a+1))\Gamma(\tfrac{1}{2}(f+a))} + \frac{\Gamma(\tfrac{1}{2}(f-a+1))\Gamma(\tfrac{1}{2}(e-a))}{\Gamma(\tfrac{1}{2}(f+a+1))\Gamma(\tfrac{1}{2}(e+a))}
\right) \,.
\end{align*}
After rearrangement in \eqref{eq:lambda'=1}, this gives the result leading to \eqref{eq:B11mumununu}:
\begin{align*}
\Bc(1,1,\mu,\mu,\nu,\nu) &= \frac{1}{2} - \frac{\Gamma(\mu+\tfrac{1}{2})\Gamma(\nu+\tfrac{1}{2})\Gamma(\mu+\nu+1)}{4\sqrt{\pi}\Gamma(1+\mu)\Gamma(1+\nu)\Gamma(\mu+\nu+\tfrac{1}{2})} \\
&= \frac{1}{2} -\frac{(\mu+\nu)B(\mu+\nu,\mu+\nu+1)}{4\mu\nu B(\mu,\mu+1)B(\nu,\nu+1)}
\,.
\end{align*} 

\noindent\textbf{Case 2:} $\lambda=\lambda'=1$, $\mu'=\mu+1$, $\nu'=\nu+1$.\\
Similarly to Case 1, this can be evaluated using the $i=0$, $j=-1$ case of~\cite{Lavoie:1996}:
i.e., for $a+b=0$, $e+f=2c+1$
\begin{align*}
\pFq 32{a, b, c}{e, f}{1} &= \frac{\Gamma(e)\Gamma(f)}{2^{2a+1}\Gamma(e-a)\Gamma(f-a)} \\
&\qquad \times\left(
\frac{\Gamma(\tfrac{1}{2}(e-a))\Gamma(\tfrac{1}{2}(f-a))}{\Gamma(\tfrac{1}{2}(e+a))\Gamma(\tfrac{1}{2}(f+a))} + \frac{\Gamma(\tfrac{1}{2}(e-a+1))\Gamma(\tfrac{1}{2}(f-a+1))}{\Gamma(\tfrac{1}{2}(e+a+1))\Gamma(\tfrac{1}{2}(f+a+1))}
\right) \,.
\end{align*}
After rearrangement, this gives
\begin{align*}
\Bc(1,1,\mu,\mu+1,\nu,\nu+1) &=\frac{4\mu\nu+3(\mu+\nu)+2}{2(2\mu+1)(2\nu+1)}
- \frac{(\mu+\nu)\Gamma(\mu+\nu+1)\Gamma(\mu+\tfrac{1}{2})\Gamma(\nu+\tfrac{1}{2})}{4\sqrt{\pi}\Gamma(\mu+1)\Gamma(\nu+1)\Gamma(\mu+\nu+\tfrac{3}{2})} \\
&	= \frac{4\mu\nu+3(\mu+\nu)+2}{2(2\mu+1)(2\nu+1)} -\frac{(\mu+\nu)B(\mu+\nu+1,\mu+\nu+2)}{\mu\nu B(\mu,\mu+1)B(\nu,\nu+1)}\,.
\end{align*}
which may be written for integer parameters as
\begin{equation}\label{eq:B11mm1nn1}
\Bc(1,1,m,m+1,n,n+1) =\frac{4mn+3(m+n)+2}{2(2m+1)(2n+1)}
- \frac{m+n}{4(m+n)+2}\frac{\binom{m+n}{m}^2}{\binom{2m+2n}{2m}}.
\end{equation}

\noindent\textbf{Case 3:} $\lambda=\lambda'=1$, $\nu=\mu$, $\nu'=\mu'$. \\
Here we consider the situation in which the random variables $Y_{\mu,\mu'}$ and $Z_{\nu,\nu'}$ are identically distributed. The corresponding hypergeometric functions in \eqref{eq:lambda'=1} may be computed using the $i=-1$, $j=1$ case of a contiguous generalization of Dixon's theorem presented in \cite{Lavoie:1994}: i.e., for $e=a-b$, $f=a-c+1$
\begin{align*}
\pFq 32{a, b, c}{e,f}{1} &= \frac{\Gamma(e)\Gamma(f)}{4^c\Gamma(f-c)\Gamma(f-b)} \\
&\qquad \times\left(
\frac{\Gamma(\tfrac{1}{2}(f-c))\Gamma(\tfrac{a}{2}-b-c+1)}{\Gamma(\tfrac{1}{2}(a+1))\Gamma(\tfrac{1}{2}(e-b))} + \frac{\Gamma(\tfrac{1}{2}(f-c+1))\Gamma(\tfrac{a}{2}-b-c+\tfrac{1}{2})}{\Gamma(\tfrac{a}{2})\Gamma(\tfrac{1}{2}(e-b+1))}
\right) \,.
\end{align*}
Taking $a=2\mu+1$, $b=\mu$ and $c=1-\mu'$, a little rearrangement gives
\begin{equation}\label{eq:B11mnmn}
\Bc(1,1,\mu,\mu',\mu,\mu')  = \frac{1}{\mu+\mu'}\left(\mu' - \frac{2 B(2\mu,2\mu')}{B(\mu,\mu')^2}\right) \,,
\end{equation}
and the moment result in \eqref{eq:B11munumunu} follows immediately.

\medskip
\noindent\textbf{Case 4:}  $\nu=\mu$, $\mu'=\nu'=\tfrac{1}{2}$. \\
For our final example we return to~\eqref{eq:ETFF} and note that for this particular set of parameters one may use an instance of Clausen's formula~(\citeauthor[16.12.2]{DLMF}) 
\[
\pFq 21{\tfrac{1}{2},\mu}{\mu+1}{t}^2 = \pFq 32{1,\mu+\tfrac{1}{2},2\mu}{\mu+1,2\mu+1}{t}
\]
to obtain
\begin{equation*}
\Bc(\lambda,\lambda',\mu,\tfrac{1}{2},\mu,\tfrac{1}{2}) =
\frac{B(\lambda+2\mu,\lambda')}{\mu^2 B(\lambda,\lambda')B(\mu,\tfrac{1}{2})^2}\,
\pFq 43{1, \mu+\tfrac{1}{2},2\mu,\lambda+2\mu}{\mu+1,2\mu+1,\lambda+\lambda'+2\mu}{1}\,.
\end{equation*}

\medskip
References \citep{Lavoie:1994,Lavoie:1996} contain numerous contiguous generalizations of the Dixon and Whipple theorems respectively, so it is clear that many more special cases could be derived. For instance, if $\lambda'=1$ and the parameters $i_1=\mu'-\mu-\lambda$, $i_2=\nu'-\nu-\lambda$ and $j=\mu+\nu-(\mu'+\nu')+\lambda$  are integers in the interval $[-3,3]$, the hypergeometric functions in~\eqref{eq:lambda'=1} can be evaluated using the $(i_1,j)$ and $(i_2,j)$ generalizations of Whipple's theorem obtained in~\cite{Lavoie:1996}. Details are left to the reader.

\section{Applications}\label{sec:path_counting}

In this final section we restrict attention to the case when all six parameters in our integral \eqref{eq:B_int_defn} are positive integers, and show that the integrals are related to some attractive problems involving a random walk and string enumeration. Let us now write $\Bc(\ell,\ell',m,m',n,n')$ to emphasise this restriction.

It is well known that the order statistics of the uniform distribution on $[0,1]$ follow beta distributions: the $m$'th smallest of $m+m'-1$ independent $\Uniform[0,1]$ random variables is distributed as $\Beta(m,m')$. As in Section~\ref{sec:intro}, let $X_{\ell,\ell'}\sim\Beta(\ell,\ell')$, $Y_{m,m'}\sim \Beta(m,m')$ and $Z_{n,n}\sim \Beta(n,n')$ be independent. Then we see that the probability $\prob{X_{\ell,\ell'}>\max\{Y_{m,m'},Z_{n,n'}\}}$ concerns $L+M+N$
i.i.d.\ uniform draws from $[0,1]$ taken in three batches, which we will label `lime', `magenta' and `navy', of sizes $L=\ell+\ell'-1$, $M=m+m'-1$ and $N=n+n'-1$ respectively, and each ordered from smallest to largest. Our interest is in the event that the $\ell$'th lime value should exceed the maximum of the $m$'th magenta and $n$'th navy values. As the individual draws from $[0,1]$ are i.i.d., every reordering of the $L+M+N$ values obtained is equally probable. Therefore (recalling \eqref{eq:prob_interpretation})

\begin{equation}\label{eq:fromBtoT}
\Bc(\ell,\ell',m,m',n,n') = \binom{L+M+N}{L,\, M,\, N}^{-1}\,T(\ell,\ell',m,m',n,n'),
\end{equation}
where $T(\ell,\ell',m,m',n,n')$ is the number of permutations of $L$ lime, $M$ magenta and $N$ navy counters (distinguishable only by their colour) placed in a line, in which the $\ell$'th lime counter lies to the right of at least $m$ magenta and $n$ navy counters. We now present our two applications.

\subsection{Exit time of a random walk on the lattice}\label{ssec:random_walk}

Consider the case where $\ell = \ell' =1$ and we only have one lime counter. Conditional upon the order of the magenta and navy counters, the lime counter is uniformly distributed on the set of `gaps' between them. The arrangements of the magenta and navy counters are in bijection with paths on a $M\times N$ square lattice, starting from position $(0,0)$ and terminating at $(M,N)$, in which single steps may only be taken to the East (magenta counter) or North (navy counter). Coordinates are specified with the eastward position first. These can be thought of as paths of a random walk on the lattice, started from the origin, which at each step moves either East or North with equal probability and which is conditioned to pass through the point $(M,N)$.

Let $A=(A_0,A_1,\dots,A_{M+N})$ be a path chosen uniformly at random from the set of all such paths, and let $U$ be uniformly distributed on the set $\{0,1,\dots,M+N\}$ (independently of $A$). Recalling equation~\eqref{eq:min_not_max}, $\Bc(1,1,m,m',n,n')$ is the probability that the lime counter lies to the left of at least $m$ magenta and $n$ navy counters. Let $R$ denote the rectangle with corners $(0,0)$ and $(m'-1,n'-1)$; once a path has left $R$ it can never re-enter, because only northward or eastward steps are permitted. Thanks to the correspondence outlined above, we see that 
\[ 
\Bc(1,1,m,m',n,n') = \prob{A_U \in R} = \prob{T_R > U} \,,
\] 
where $T_R = \min\{k: A_k \notin R\}$ is the time at which the random walk exits $R$. Since $U$ is uniform this may equivalently be written as 
\[
\Bc(1,1,m,m',n,n') = \frac{1}{M+N+1} \sum_{k=0}^{M+N}\prob{T_R > k} = \frac{\mathbb{E}[T_R]}{M+N+1} \,,
\]
and so we see that our formulae for $\Bc(1,1,m,m',n,n')$ may be used to calculate the expected exit time of this conditioned random walk.

For example, in the case $m=n = m'-1=n'-1$ we may obtain the asymptotic behaviour as $n\to\infty$ of the exit time from the set $[0,n]^2$ of a random walk conditioned to pass through the point $(2n,2n)$:
\begin{align*}
\mathbb{E}[T_{[0,n]^2}] &= (4n-1)\Bc(1,1,n,n+1,n,n+1) \\ 
&= (4n-1)\left(
\frac{4n^2+6n+2}{2(2n+1)^2}
- \frac{n}{4n+1}\frac{\binom{2n}{n}^2}{\binom{4n}{2n}}\right) = 2n - \sqrt{\frac{2n}{\pi}} + O(1)\,,
\end{align*}
thanks to \eqref{eq:B11mm1nn1} and an application of Stirling's formula. This can be compared to a simple random walk on the lattice which at each step moves either East or North but which is \emph{not} conditioned to pass through the point $(2n,2n)$. The corresponding expected exit time is given by:
\[
\sum_{k=n+1}^{2n+1} 2k\, \frac{\binom{k-1}{n}}{2^k} = 2(n+1) - \frac{2 \Gamma(n+3/2)}{\sqrt{\pi}\Gamma(n+1)} = 2n - 2\sqrt{\frac{n}{\pi}} +O(1)\,,
\] 
as is easily seen by counting the number of paths which occupy a boundary point of $R$ at time $k-1$ and subsequently exit at time $k$.

\subsection{String enumeration}\label{ssec:strings}

It follows from the discussion leading to \eqref{eq:fromBtoT} that formulae for $\Bc(\ell,\ell',m,m',n,n')$ yield solutions to a natural combinatorial problem, which may also be expressed as the enumeration of strings of certain types. Specifically, $T(\ell,\ell',m,m',n,n')$
is the number of strings comprised of $L$ {\tt X}'s, $M$ {\tt Y}'s and $N$ {\tt Z}'s, in which the $\ell$'th {\tt X} appears after at least $m$ {\tt Y}'s and $n$ {\tt Z}'s. 

To conclude, we give closed formulae for various cases of this problem, drawing upon some of the special cases considered in Section~\ref{ssec:special}. For example, $T(1,1,m,m, n,n)$ is the number of distinct strings consisting of $2m-1$ {\tt Y}'s, $2n-1$ {\tt Z}'s and one {\tt X} in which the character {\tt X} lies to the right of more than half of the {\tt Y}'s and more than half of the {\tt Z}'s; using \eqref{eq:B11mmnn} one has
\begin{align*}
T(1,1,m,m,n,n) &=  \frac{(2m+2n-1)!}{(2m-1)!(2n-1)!}\left(
\frac12 - \frac14 \frac{\binom{m+n}{m}^2}{ \binom{2m+2n}{2m} } 
\right)\\
& = \frac{1}{2B(2m,2n)} - \frac{1}{2B(m,n)}\binom{m+n}{m}\,.
\end{align*}
We tabulate some values in Table~\ref{tab:T11mmnn}. The left-most column forms the sequence 
$T(1,1,1,1,n,n)= \tfrac{1}{2}n(3n+1)$ of second pentagonal numbers, $\oeis{A005449}(n)$ in the On-Line Encyclopedia of Integer Sequences (OEIS), but the full triangle is not presently known to the OEIS, nor is the subsequence
\begin{align*}
T(1,1,n,n,n,n)  &= \frac{1}{2B(2n,2n)} - \frac{1}{2B(n,n)}\binom{2n}{n} \\
&= 2,~52,~1086,~20840,~382510,\ldots \quad (n=1,2,\ldots)
\end{align*}
obtained when there is an equal number ($2n-1$) of {\tt Y}'s and {\tt Z}'s. 
\begin{table}\begin{center}
		\begin{tabular}{|c|ccccc|}\hline
			& \multicolumn{5}{c|}{$m$}\\ \hline 
			$m+n$ & $1$ & $2$ & $3$ & $4$ & $5$ \\
			\hline
			$2$ & $2$ & & & & \\	
			$3$ & $7$ & $7$ & & & \\ 
			$4$ & $15$ &  $52$ & $15$ && \\
			$5$ & $26$ & $192$ & $192$ & $26$ & \\
			$6$ & $40$ & $510$ & $1086$ & $510$ & $40$\\
			\hline
		\end{tabular}
	\end{center}
	\caption{Values of $T(1,1,m,m,n,n)$ for small integer values.}\label{tab:T11mmnn}
\end{table}

Similarly, using~\eqref{eq:B11mm1nn1} and rearranging, one has
\[
T(1,1,m,m+1,n,n+1)=\frac{4mn+3(m+n)+2}{4(m+n+1)} \binom{2(m+n+1)}{2m+1} - \frac{m+n}{2}\binom{m+n}{m}^2
\]
which counts the number of distinct strings consisting of $2m$ {\tt Y}'s, $2n$ {\tt Z}'s and one {\tt X} in which the character {\tt X} lies to the right of at least half of the {\tt Y}'s and at least half of the {\tt Z}'s,
tabulated in Table~\ref{tab:T11mm1nn1}. We also have
\begin{align*}
T(1,1,n,n+1,n,n+1) & =\frac{n+1}{2}\binom{4n+2}{2n+1} - n\binom{2n}{n}^2
\\ &= 
16,~306,~5664,~101950,~1798776,\ldots \qquad (n=1,2,\ldots)
\end{align*}
for the case where there is an equal number ($2n$) of {\tt Y}'s and {\tt Z}'s. 

\begin{table}\begin{center}
		\begin{tabular}{|c|ccccc|}\hline
			& \multicolumn{5}{c|}{$m$}\\ \hline 
			$m+n$ & $1$ & $2$ & $3$ & $4$ & $5$ \\
			\hline
			$2$ & $16$ & & & & \\
			$3$ & $53$ & $53$ & & & \\	
			$4$ & $124$ & $306$ & $124$ & & \\ 
			$5$ & $240$ & $1103$ &  $1103$ & $240$ & \\
			$6$ & $412$ & $3043$ & $5664$ &  $3043$ &  $412$ \\ 
			\hline
		\end{tabular}
	\end{center}
	\caption{Values of $T(1,1,m,m+1,n,n+1)$ for small integer values.}\label{tab:T11mm1nn1}
\end{table}

For a final example we use \eqref{eq:B11mnmn} to count the number of distinct strings of $m+n-1$ {\tt Y}'s, $m+n-1$ {\tt Z}'s and one {\tt X} in which the character {\tt X} lies to the right of at least $m$ of the {\tt Y}'s and $m$ of the {\tt Z}'s: 
\begin{align*}
T(1,1,m,n,m,n) &= \frac{(2m+2n-1)!}{(m+n-1)!^2}\frac{1}{m+n}\left(n - \frac{2 B(2m,2n)}{B(m,n)^2}\right) \\
&= n \binom{2 m + 2 n - 1}{m + n} - \frac{2 m n}{m+n} \binom{2 m - 1}{m}\binom{2 n - 1}{n} \,. 
\end{align*}
We table some values of $T(1,1,m,n,m,n)$ in Table~\ref{tab:T11mnmn}. The left-most column forms the sequence 
$T(1,1,1,n,1,n)= 2n \binom{2n}{n-1}=\oeis{A253487}(n-1)$, and the leading diagonal forms the sequence
$T(1,1,m,1,m,1)=  \binom{2m}{m}= \oeis{A000984}(m)$ of central binomial coefficients, but once again the full triangle is not presently known to the OEIS.

\begin{table}[h]\begin{center}
		\begin{tabular}{|c|ccccc|}\hline
			& \multicolumn{5}{c|}{$m$}\\ \hline 
			$m+n$ & $1$ & $2$ & $3$ & $4$ & $5$ \\
			\hline
			$2$ & $2$ & & & & \\
			$3$ & $16$ & $6$ & & & \\	
			$4$ & $90$ & $52$ & $20$ & & \\ 
			$5$ & $448$ & $306$ &  $180$ & $70$ & \\
			$6$ & $2100$ & $1568$ & $1086$ &  $644$ &  $252$ \\ 
			\hline
		\end{tabular}
	\end{center}
	\caption{Values of $T(1,1,m,n,m,n)$ for small integer values.}\label{tab:T11mnmn}
\end{table}

These special cases have been picked out because the final answers are particularly simple. However, any of the formulae given for $\Bc(\ell,\ell',m,m',n,n')$ provide an answer to a corresponding 
combinatorial problem; other instances are left to the reader.

\bibliographystyle{chicago}
\bibliography{beta}

\end{document}